\documentclass[11pt]{article}
\usepackage{graphicx,amsmath,amssymb,mathrsfs,tcolorbox}
\usepackage[a4paper,  margin= 1 in]{geometry}
\usepackage[utf8]{inputenc}
\usepackage[english]{babel}
\usepackage[utf8]{inputenc}
\usepackage{amsmath}
\usepackage{float}
\usepackage{amsfonts}
\usepackage{amssymb} 
\usepackage{authblk}
\usepackage{natbib}
\usepackage[none]{hyphenat}
\usepackage{setspace}
\setcitestyle{authoryear,open={(},close={)}}

\newtheorem{theorem}{Theorem}[section]

\newtheorem{lemma}[theorem]{Lemma}

\begin{document}

\title{\textbf{ Correction factor of FWER for normal distribution in nearly independent setup }}

\author[1]{Nabaneet Das}
\author[2]{Subir K. Bhandari}
\affil[1]{Indian Statistical Institute,Kolkata}
\affil[2]{Indian Statistical Institute,Kolkata}

\maketitle


\begin{abstract} 
In this paper, we have attempted to study the behaviour of the family wise error rate (FWER) for Bonferroni's procedure in a nearly independent setup for normal distribution. In search for a suitable correlation penalty, it has been noted that the root mean square (RMS) of correlations is not appropriate under this setup as opposed to the study of \cite{efron2007correlation}. We have provided a suitable correction factor for deviation from independence and approximated the FWER under this nearly independent setup. 

\end{abstract}

\section{Introduction}
\label{intro}

Dependence among hypotheses in simultaneous testing problem have caused a great concern among researchers. Although many efforts have been made to generalize the existing methods under dependence (\cite{yekutieli1999resampling},
\cite{benjamini2001control}, 
\cite{sarkar2002some}, \cite{sarkar2008methods},
\cite{efron2007correlation}, \cite{efron2012large} etc.), very few literature is available which explicates the effect of dependence on the existing methods. \cite{sarkar2008methods} discusses false discovery rate (FDR) control under dependence. \cite{efron2007correlation} in his study of empirical Bayes methods, has shown that, the correlation penalty depends on the root mean square (RMS) of correlations. An excellent review of the whole literature can be found in \cite{efron2012large}. \\ 
When correlation is present, these methods show some undesirable characteristics (Being too conservative when positive correlation is present in some cases). (see \cite{yekutieli1999resampling}, \cite{sarkar2008methods}, \cite{genovese2006false}, \cite{sun2009large}) The conservative nature of these methods result in loss of power. Hence, it is imperative to study the error rate of these methods in order to rectify the effect of dependence. In a recent paper \cite{das2021bound} have shown that under equicorrelated normal setup with correlation $\rho$, the widely used Bonferroni method asymptotically controls the family-wise error rate (FWER) at level $ \alpha (1 - \rho)$ instead of $ \alpha$ and \cite{dey2021behaviour} further improves this upper bound by showing that, the asymptotic FWER is zero if correlations are bounded away from zero. \cite{dey2021behaviour} also studies the FWER under normal distribution in more general setups (non-asymptotic and non-equicorrelated setups) and confirm the conservative nature of Bonferroni's method. In a simulation study, \cite{das2020observation} have pointed out the exact behaviour of FWER and FDR of Bonferroni's method and Benjamini-Hochberg FDR control for moderate to large number of hypotheses under normal distribution. \\ 
In this work, we have also considered jointly normal distribution and we shall focus on approximating the FWER of Bonferroni's method under "nearly independent" setup. We have provided an asymptotic correction factor to the FWER of Bonferroni's method under this setup.\\ 
This article is organized as follows : in section 2, we introduce the framework of our work. In section 3, we have shown the motivation of considering the "nearly independent" setup and in section 4 we have approximated the FWER under this setup and provided a correction factor which accounts for the deviation from independence. 

\section{Description of the problem} 
Let $ X_1,X_2,...... $ be a sequence of observations and the null hypotheses are  
 
$$H_{0i} : X_i \sim N(0,1) \: \: i=1,2,.... $$  
Here we have considered one sided tests (This means, $H_{0i}$ is rejected for large values of $X_i$ (say $X_i > c$)). A classical measure of the type-I error is FWER, which is the probability of falsely rejecting at least one null hypothesis (Which happens if $ X_i > c $ for some $i$ and the probability is computed under the intersection null hypothesis $ H_0= \bigcap\limits_{i=1}^{n} H_{0i} = \{ X_i \sim N(0,1) \: \: \forall \: i=1,2,...,n \}$).
Then,
\begin{center} 
\textbf{FWER} = P(At least one false rejection)= P($ \bigcup\limits_{i=1}^{n} \{ X_i > c \} $ $| $  $ H_0 $ )  
\end{center}
Suppose, $ Corr(X_i,X_j) = \rho_{ij} \: \: \:  \forall \: i \neq j $ and $\Sigma$ is the corresponding correlation matrix. \\ 


We assume the following about the correlation matrix $ \Sigma $. 

$$\rho_{ij} = O \: ( \frac{1}{n^{ \beta } }  )  \: \: \forall \: i \neq j \text{ for some } \beta > 0 $$ 
For convenience, we shall call this setup as "\textbf{nearly independent setup}".\\
\vspace{0.03 in} 
Suppose $ P (X_i > c  \: | \:  H_{oi} ) = \alpha_n$. Then, under the assumption of independence of hypotheses (i.e. if $\Sigma $ is identity matrix), we have $ FWER  = 1 -  ( 1 - \alpha_n)^n $. We shall study the value of the FWER under nearly independent setup. Bonferroni's method is based on $ \alpha_n = \frac{ \alpha}{n}$. So, we shall consider the setup where $\lim_{ n \to \infty} n \alpha_n = \alpha$. 

\section{Motivation behind the nearly independent setup} 
\cite{efron2007correlation} in his study of empirical Bayes methods, has pointed out that the correlation penalty on the summary statistics based on empirical c.d.f. depends on the root mean square (RMS) of correlations ($\frac{1}{n(n-1)} \sum\limits_{i \neq j } \sum \rho_{ij}^2 $). A natural question might arise from this result is that "does the RMS of correlations act as a correlation penalty" in our framework also ?" \\ 
If this were true, then we would have $\text{ FWER } \sim \text{ FWER under independence} $ as $\frac{1}{n(n-1)} \sum\limits_{i \neq j } \sum \rho_{ij}^2  \to 0$. \\ 
We can construct an easy counterexample to answer this question. \\ 
Let $M_n ( \rho) $ denote the $n*n $ equicorrelation matrix with correlation $\rho$ and $\mathbf{0}_n$ denote the $n*n$ matrix with all zero entries. Consider the $n^2*n^2$ block-diagonal matrix with $n$ many $M_n ( \rho) $ blocks in the diagonal. 
$$ 
\text{ i.e.  } \Sigma = 
\begin{pmatrix} 
M_n ( \rho) & \mathbf{0}_n & ... & \mathbf{0}_n \\
\mathbf{0}_n & M_n ( \rho) & ... & \mathbf{0}_n \\
...& ..... & & ..... \\ 
\mathbf{0}_n & \mathbf{0}_n & ... & M_n ( \rho)
\end{pmatrix}
$$
This is the correlation matrix of the $n^2$ variables $(X_1, ..., X_{n^2} ) $ such that, 
\begin{itemize} 
\item $ (X_1, .., X_n) , (X_{n+1} , .., X_{2n}) , ...., (X_{n^2 - n +1 } , ..., X_{n^2} ) $ are independent. 
\item In each block $ (X_{(k-1)n +  1 } , .., X_{kn})$  we have $ Corr (X_i, X_j ) = \rho $ for $i \neq j$.  (for $k=1, \dots, n $)

\end{itemize} 
If we apply Bonferroni ($\alpha$) method in this setup then $ \alpha_n = \frac{\alpha}{n^2}$. By the similar approach of \cite{das2021bound}, we can argue that,\\ 
$ 1 - FWER \geq ( (1 - \frac{ \alpha}{n^2}) -  ( 1 -\rho) [(1 - \frac{\alpha}{n^2}) - (1 - \frac{\alpha}{n^2})^n])^n  = ( 1 - \frac{ \alpha}{n^2})^n [ 1 - (1 - \rho) [1 - ( 1 - \frac{ \alpha}{n^2} )^{n-1}]]^n $  \vspace{0.03 in} 
Since $ 1 - (1-x)^k \leq kx \: \: \: \forall \: 0 < x < 1 $, this implies, \\   $$  
1 - FWER  \geq ( 1- \frac{ \alpha}{n^2})^n [ 1 - \frac{ (1- \rho)(n-1) \alpha}{ n^2} ]^n \to e^{ - \alpha (1 - \rho)}  \text{ as }  n \to \infty  $$. \\ 
Under this setup, we finally have $ FWER  \leq \alpha ( 1 - \rho)$ asymptotically. It is interesting to note that, 
\begin{itemize}
    \item In this setup, mean of the absolute values of correlations is $ \frac{ n \binom{n}{2} \rho}{ \binom{n^2}{2} } = O( \frac{1}{n} )   $ and so is the mean square of correlations. Although the root mean square (RMS) of correlations becomes smaller and smaller, the FWER is not close to the FWER under independence. So, the correlation penalty is neither dependent on RMS of correlations nor the mean of their absolute values.  
    \item It is also interesting to note that, whenever $\max\limits_{ i \neq j } | \rho_{ij} | $ is bounded away from 0, it is possible to find a setup where FWER widely differs from the one obtained under independence. 
\end{itemize} 
So, if we ask ourselves "how far can we go from independence without losing much?", then the above example suggests us to bound the absolute values of the correlations. Under the nearly independent setup, we bound the absolute values by the order $ \frac{1}{n^{ \beta } } $ for some $ \beta > 0 $.  

\section{FWER under the nearly independent setup} 
Let's denote FWER under correlation matrix $ \Sigma $ by $ FWER ( \Sigma) $.
\begin{theorem}{}\label{4.1}
Fix a sufficiently large integer $K$. Then, under the nearly independent setup, FWER can be approximated by the following way. 
$$ FWER ( \Sigma ) 
  \sim \sum\limits_{i=1}^{ K} \frac{ (-1)^{ i-1 } \alpha^i}{i! } +  
 \frac{c^2 \bar{ \rho} }{2} \sum\limits_{ i=2}^{K } (-1)^{i-1} \frac{ \alpha^i}{ (i-2) ! }   $$  
 Where $ \bar{ \rho} = \frac{1}{n (n-1) } \sum\limits_{ i \neq j } \sum \rho_{ij} $ is the mean of correlations. 
\end{theorem} 
\textbf{\underline{Remarks}} :- 
 It is interesting to note that, $ \sum\limits_{i=1}^{ K} \frac{ (-1)^{ i-1 } \alpha^i}{i! }  \approx 1 - e^{ - \alpha } $  for sufficiently large $K$. Also $ \lim_{ n \to \infty} 1 - ( 1- \alpha_n)^n = 1 - e^{ - \alpha }  $ is the limiting form of FWER under independence. So, if we use large enough $K$ for this approximation then the first term is nearly equal to the FWER for independent case. So, the second term $\frac{c^2 \bar{ \rho} }{2} \sum\limits_{ i=2}^{K } (-1)^{i-1} \frac{ \alpha^i}{ (i-2) ! } $ acts as a correction factor or the amount of deviation from independence. Under the nearly independent setup, $ \bar{ \rho} = O ( \frac{1}{ n^{ \beta } } )$ (The correction factor is very small) 
 \\
 \\
Before we proceed with the proof of theorem 4.1, the following important theorem must be stated. 
 \begin{tcolorbox}

 \begin{theorem}{ \underline{Multivariate Mill's Ratio (\cite{savage1962mills})} }

Let $\underset{\sim}{X} \sim N_m ( \mathbf{0}_m , V) $ and $ M = V^{-1}$.  \\ 
$ F( \underset{\sim}{a} , M ) =  P( \underset{\sim}{X} > \underset{\sim}{a} )  $ and $f( \underset{\sim}{a} , M ) = \frac{|M|^{ \frac{1}{2} } }{ (2 \pi)^{ \frac{k}{2} } } \exp ( - \frac{1}{2} \underset{\sim}{a}^T M \underset{\sim}{a} )  $ \\ 

Let $ \underset{\sim}{\Delta} = \underset{\sim}{a}^T M $ (i.e. $ \Delta_i = \sum\limits_{j} a_j m_{ji} \: \: \: \:  i=1,2,..., m $)  \\ 

If $ \Delta_i > 0 \: \: \forall \: i , $ then  

$$ 1 - \frac{1}{2} \sum\limits_{i} \sum\limits_{j} \frac{ m_{ij} ( 1 + \delta_{ij} ) }{ \Delta_i \Delta_j } < \frac{ F( \underset{\sim}{a} , M )}{f( \underset{\sim}{a} , M ) (\prod\limits_{i=1}^{m} \Delta_i )^{-1}}  < 1 $$  

Here $ \delta_{ij} = $ Kronecker's Delta 
\end{theorem}

\end{tcolorbox}

\textbf{\underline{Proof of theorem 4.1}} :- \\

$
\mathbf{FWER} \:  ( \Sigma) 
 = \sum\limits_{i=1}^{n} P(X_i > c ) - \sum\limits_{ i \neq j } \sum P(X_i > c , X_j > c )  
 + \sum \sum\limits_{ i,j,k \text{ distinct}}  \sum P(X_i > c , X_j > c , X_k > c ) .... + (-1)^{n-1} P(X_1 > c , ..., X_n > c ) $

Fix any $ 1 \leq k \leq K < < n $ where $K $ is a fixed positive integer. 
\\
\begin{lemma}{} :- 
For any $  1 \leq i_1 < ... < i_k \leq n $. Then, 
$$ P(X_{i_1 } > c , ..., X_{i_k} > c ) \sim f ( c \mathbf{1}_k ) ( \prod\limits_{i=}^{k} \Delta_i)^{-1}  $$ \\ 
(Where $ \Delta = c \mathbf{1}_k^{T} W^{-1} , \: \: W $ is the correlation matrix of $ (X_{i_1} , X_{i_2} , ..., X_{ik} )$ and $f(.)$ is the density function of $N_k ( \mathbf{0}_k , W) $ distribution)  
\end{lemma} \vspace{0.02 in} 
\underline{Proof } :-   \\ 
Let $W = I + R $. \\ 
$ \Rightarrow W^{-1 } = ( I + R)^{-1} = I - R + o (R)  $\\
By the nearly independent assumption, we have $ ( I +R )^{ -1 } \approx ( I - R ) $. \\ 
So, $ \Delta = c \mathbf{1}_k^{T} W^{-1} \approx c \mathbf{1}_k^{T}  ( I - R)  $. \\ 
$ \Rightarrow \Delta_i = c ( 1 - \sum\limits_{ j \neq i } \rho_{ji} )  > 0   \: \: \: \: \forall \: i $ (By the nearly independent assumption)  \\ 
Now it is clear that, $ 1 - \frac{1}{2} \sum\limits_{i} \sum\limits_{j} \frac{ m_{ij} ( 1 + \delta_{ij} ) }{ \Delta_i \Delta_j }  = 1 - O ( \frac{1}{c^2} )$. (This is because the sum is over a fixed no. of terms ($k^2$ many terms)). \\ 
Since $ c \to \infty $ as $ n \to \infty $, we can say that, 
$  \lim_{ n \to \infty } \frac{P(X_{i_1 } > c , ..., X_{i_k} > c ) }{ f ( c \mathbf{1}_k ) ( \prod\limits_{i=}^{k} \Delta_i)^{-1}} = 1 $ and hence the result.

\underline{ \textbf{Approximation of } $f ( c \mathbf{1}_k ) ( \prod\limits_{i=}^{k} \Delta_i)^{-1}  $ } \\

\vspace{0.05 in} 
By the nearly independent assumption on $ \Sigma$, we can say that, $ | I +R | \sim 1 $.\\ 

Observe that,  

\vspace{0.02 in}
$ f ( c \mathbf{1}_k) \sim \frac{1}{ (2 \pi)^{ \frac{k}{2} } } \exp ( - \frac{c^2}{2} \mathbf{1}_k^T (I - R) \mathbf{1}_k) = 
\frac{1}{ (2 \pi)^{ \frac{k}{2} } } \exp ( - \frac{kc^2}{2} ) \exp ( \frac{c^2}{2} \sum\limits_{ l \neq m } \sum \rho_{ lm} )  $ \\ 

\vspace{0.02 in}
Since $ \rho_{lm} =  O( \frac{1}{n^{ \beta} } )  $ and $c^2 = O ( \log n) $, this implies, $ c^2 \sum\limits_{ l \neq m } \sum \rho_{lm} = o(1) $. \\

We know that, $ e^x \sim 1 + x $ for sufficiently small $x$. So, 
$$ \exp ( \frac{c^2}{2}\sum\limits_{ l \neq m } \sum \rho_{lm} ) \sim 1 + \frac{c^2}{2}\sum\limits_{ l \neq m } \sum \rho_{lm} $$   

$ \Rightarrow 
f( c \mathbf{1}_k) \sim  
\frac{1}{ (2 \pi)^{ \frac{k}{2} } } \exp ( - \frac{kc^2}{2} ) ( 1 + \frac{c^2}{2}\sum\limits_{ l \neq m } \sum \rho_{lm}  ) $   and $ ( \prod\limits_{i=1}^{k} \Delta_i) = c^k \prod\limits_{i=1}^{k} (1 - \sum\limits_{ j \neq i } \rho_{ji} ) $. \\ 

Since $ \rho_{ij} = O ( \frac{1}{n^{ \beta } } ) $, this implies, $ \prod\limits_{i=1}^{k} (1 - \sum\limits_{ j \neq i } \rho_{ji} )  \sim 1 $.  \\ 

Hence,  
$$f ( c \mathbf{1}_k ) ( \prod\limits_{i=}^{k} \Delta_i)^{-1}   \approx  \frac{c^k (  1 + \frac{c^2}{2}\sum\limits_{ l \neq m } \sum \rho_{lm} ) }{ ( 2 \pi)^{ \frac{k}{2} } \exp ( \frac{kc^2}{2} ) }$$  
 
Recall that, $ \Phi ( -c) = \alpha_n \sim \frac{ \alpha}{n} $.
\\

If $ \phi(.) $ denote the standard normal density, then for large enough $c$, then 
$$ \phi (c) \sim c \Phi (-c)   \Rightarrow n \sim \alpha \sqrt{ 2 \pi }  c e^{ \frac{c^2}{2} } $$   

This implies, for any finite $k$, 

\vspace{0.03 in} 

$ f ( c \mathbf{1}_k ) ( \prod\limits_{i=}^{k} \Delta_i)^{-1} \sim  ( \frac{ \alpha}{n} )^k ( 1 + \frac{c^2}{2} \sum\limits_{ l \neq m } \sum \rho_{ lm} ) $  

In particular, 
$$ P( X_{i_1}  > c , ..., X_{ik } > c )   \sim 
( \frac{ \alpha}{n} )^k ( 1 + \frac{c^2}{2} \sum\limits_{ l \neq m  \in \{i_1, .., i_k \}} \sum \rho_{ lm} )
$$ 

\begin{tcolorbox}
Approximation of $ FWER ( \Sigma ) $ under the nearly independent setup 
\end{tcolorbox}

$
\mathbf{FWER} \:  ( \Sigma) 
 = \sum\limits_{i=1}^{n} P(X_i > c ) - \sum\limits_{ i \neq j } \sum P(X_i > c , X_j > c )  
 + \sum \sum\limits_{ i,j,k \text{ distinct}}  \sum P(X_i > c , X_j > c , X_k > c ) .... + (-1)^{n-1} P(X_1 > c , ..., X_n > c ) $  \\ 
 
 Clearly, $ P(X_i > c) = \Phi ( -c) = \alpha_n $ and hence  $ \sum\limits_{i=1}^{n} P(X_i > c ) = n \alpha_n \sim \alpha $ 
 
 For a fixed positive integer $K $ ($ K > 2 $), we have  \\ 
 $ FWER ( \Sigma ) = n \alpha_n + \sum\limits_{ i=2}^{K} \sum\limits_{ \substack{ u_1, ..., u_i \\  \text{ distinct} } }  ( -1)^{ i-1} P( X_{u_1} > c , ..., X_{u_i} > c) + \sum\limits_{i=K+1}^{n} \sum\limits_{ \substack{ u_1, ..., u_i  \\ \text{ distinct } } } (-1)^{ i-1} P( X_{u_1} > c , ..., X_{u_i } > c )  $   \\ 
 
 Clearly,

$\sum\limits_{ i=2}^{K} \sum\limits_{ u_1, ..., u_i \text{ distinct} } ( -1)^{ i-1} P( X_{u_1} > c , ..., X_{u_i} > c)  $  \\ 
 
$ \sim \sum\limits_{ i=2}^{K} \sum\limits_{ u_1, ..., u_i \text{ distinct} } ( -1)^{ i-1} ( \frac{ \alpha}{n})^i ( 1 + \frac{c^2}{2} \sum\limits_{ l \neq m  \in \{ u_1, ..., u_i  \}} \sum \rho_{ lm} )  $\\ 
 
$ \sim \alpha +   \sum\limits_{i=2}^{K} (-1)^{ i-1} \frac{ \binom{n}{i} \alpha^i}{n^i }  + \frac{c^2}{2} \sum\limits_{i=2}^{K} (-1)^{i-1} ( \frac{ \alpha}{n})^i \binom{n-2}{i-2} (\sum\limits_{ i \neq j } \sum \rho_{ij} ) $\\ 
 
$ \sim \sum\limits_{i=1}^{ K} \frac{ (-1)^{ i-1 } \alpha^i}{i! } +  
 \frac{c^2}{2 n^2 } \sum\limits_{ i=2}^{K } (-1)^{i-1} \frac{ \alpha^i}{ (i-2) ! }  (\sum\limits_{ i \neq j } \sum \rho_{ij} ) $ \\ 
 
$ \sim \sum\limits_{i=1}^{ K} \frac{ (-1)^{ i-1 } \alpha^i}{i! } +  
 \frac{c^2 \bar{ \rho} }{2  } \sum\limits_{ i=2}^{K } (-1)^{i-1} \frac{ \alpha^i}{ (i-2) ! }  $
 
 \vspace{0.05 in}

 We ignore the tail term $ \sum\limits_{i=K+1}^{n} \sum\limits_{ u_1, ..., u_i \text{ distinct }} (-1)^{ i-1} P( X_{u_1} > c , ..., X_{u_i } > c ) $ from $ FWER ( \Sigma ) $ and we approximate the FWER by the following formula. 
 
 $$ 
 FWER ( \Sigma ) 
  \sim \sum\limits_{i=1}^{ K} \frac{ (-1)^{ i-1 } \alpha^i}{i! } +  
 \frac{c^2 \bar{ \rho} }{2  } \sum\limits_{ i=2}^{K } (-1)^{i-1} \frac{ \alpha^i}{ (i-2) ! }  $$  
\textbf{\underline{Remarks}} :-  
It is interesting to note that, taking $ \rho_{ij}$'s $ O ( \frac{1}{ n^{ \beta} } ) $ is not necessary for the proof. It works for any order faster than $ \frac{1}{c^2} $ or $ \frac{1}{ \log n} $. 

\section{Simulation results} 
Bonferroni's procedure controls FWER at desired level under independence. We examine the FWER in the nearly independent setup and compare with the corrected value of FWER as per theorem \ref{4.1}. For these simulations we have considered $ n = 5000 $ and $ \beta = 0.4, 0.6, 0.8 , 1$. They are repeated on four levels of significance ($ \alpha $) ( namely $\alpha = 0.01, 0.05, 0.1, 0.2$). For each combination, the actual FWER is estimated based on 10,000 replications and the corrected value as per theorem \ref{4.1} is computed based on $K=15$.

\begin{table}[ht]
\centering
\caption{FWER for level $\alpha = 0.01$}
\begin{tabular}[t]{ |c|c|c|c| } 
 \hline
  $\beta$ & FWER & FWER under independence & FWER with correction\\ 
 \hline
 0.4 & 0.0121 & 0.00995 &  0.00993\\ 
 0.6& 0.011 & 0.00995 & 0.00995 \\
 0.8  & 0.0108 & 0.00995 & 0.00995 \\
 1 & 0.0089 & 0.00995 & 0.00995 \\ 
 \hline
\end{tabular}
\end{table}

\begin{table}[ht]
\centering
\caption{FWER for level $\alpha = 0.05$}
\begin{tabular}[t]{ |c|c|c|c| } 
 \hline
  $\beta$ & FWER & FWER under independence & FWER with correction\\ 
 \hline
 0.4 & 0.0543 & 0.0487 &  0.0484\\ 
 0.6& 0.0475 & 0.0487 & 0.0487 \\
 0.8  & 0.0422 & 0.0487 & 0.0488 \\
 1 & 0.0495 & 0.0487 & 0.0488  \\ 
 \hline
\end{tabular}
\end{table} 

\begin{table}[ht]
\centering
\caption{FWER for level $\alpha = 0.1$}
\begin{tabular}[t]{ |c|c|c|c| } 
 \hline
  $\beta$ & FWER & FWER under independence & FWER with correction\\ 
 \hline
 0.4 & 0.1077 & 0.0952 &  0.0934\\ 
 0.6& 0.0901 & 0.0952 & 0.095 \\
 0.8  & 0.0967 & 0.0952 & 0.0951 \\
 1 & 0.0983 & 0.0952 & 0.0952  \\ 
 \hline
\end{tabular}
\end{table} 

\begin{table}[ht]
\centering
\caption{FWER for level $\alpha = 0.2$}

\begin{tabular}[t]{ |c|c|c|c| } 
 \hline
  $\beta$ & FWER & FWER under independence & FWER with correction\\ 
 \hline
 0.4 & 0.0192 & 0.0181 &  0.0178\\ 
 0.6& 0.0182 & 0.0181 & 0.018 \\
 0.8  & 0.0176 & 0.0181 & 0.0181 \\
 1 & 0.018 & 0.0181 & 0.0181 \\ 
 \hline
\end{tabular}
\end{table}

The simulation results show that there is not much difference between the actual FWER and FWER under independence in the nearly independent setup. Since the correction factor is very small, the corrected value also remains very close to the value under independence. \cite{dey2021behaviour} suggests that, under the equicorrelated setup, the asymptotic FWER function have a discontinuity at 0 and then it is identically zero for any positive value of correlation. This study aims to investigate how far can we go from independence under the general setup without losing much. FWER remains equivalent to the independence in the nearly independent setup also. For a further study in this area, it would be interesting to provide an explicit correction factor due to correlations. 
\bibstyle{spbasic}

\bibliography{main.bib}

\end{document}